\documentclass{amsart}

\usepackage{amsmath,amssymb,mathrsfs}

\newcommand{\A}{{\mathbb A}}
\newcommand{\B}{{\mathbb B}}
\newcommand{\C}{{\mathbb C}}

\newcommand{\BB}{{\mathbb B}}
\newcommand{\BC}{{\mathbb C}}\newcommand{\BD}{{\mathbb D}}

\newcommand{\BR}{{\mathbb R}}

\newcommand{\BZ}{{\mathbb Z}}
\newcommand{\cA}{{\mathcal A}}\newcommand{\cB}{{\mathcal B}}
\newcommand{\cC}{{\mathcal C}}
\newcommand{\cE}{{\mathcal E}}\newcommand{\cF}{{\mathcal F}}
\newcommand{\cH}{{\mathcal H}}

\newcommand{\cK}{{\mathcal K}}\newcommand{\cL}{{\mathcal L}}

\newcommand{\cT}{{\mathcal T}}
\newcommand{\cU}{{\mathcal U}}\newcommand{\cV}{{\mathcal V}}
\newcommand{\cW}{{\mathcal W}}\newcommand{\cX}{{\mathcal X}}
\newcommand{\cY}{{\mathcal Y}}

\newcommand{\al}{\alpha}

\newcommand{\ga}{\gamma}

\newcommand{\ze}{\zeta}

\newcommand{\la}{\lambda}

\newcommand{\si}{\sigma}\newcommand{\Si}{\Sigma}

\newcommand{\vph}{\varphi}
\newcommand{\om}{\omega}\newcommand{\Om}{\Omega}

\newcommand{\diag}{\textup{diag\,}}

\newcommand{\mat}[2]{\ensuremath{\left[\begin{array}{#1}
#2
\end{array} \right]}}
\newcommand{\ov}[1]{{\overline{#1}}}

\newcommand{\inn}[2]{\ensuremath{\langle #1,#2 \rangle}}
\newcommand{\tu}[1]{\textup{#1}}

\newcommand{\half}{\frac{1}{2}}
\newcommand{\ands}{\quad\mbox{and}\quad}

\newcommand{\Sys}{\{A,B,C,D,\pi\}}
\newcommand{\col}{\text{col}}
\newcommand{\Abel}{{\bf a}}

\newtheorem{theorem}{Theorem}[section]

\newtheorem{proposition}[theorem]{Proposition}

\newtheorem{example}[theorem]{Example}

\begin{document}

\title{A $W^*$-correspondence approach to multi-dimensional linear dissipative systems}

\author{J.A. Ball}
\address{Department of Mathematics\\
Virginia Tech\\
Blacksburg VA 24061, USA} \email{ball@math.vt.edu}

\author{S. ter Horst}
\address{Department of Mathematics\\
Virginia Tech\\
Blacksburg VA 24061, USA} \email{terhorst@math.vt.edu}

\maketitle

\begin{abstract}
Recent work of the operator algebraists P.~Muhly and B.~Solel,
primarily motivated by the theory of operator algebras and
mathematical physics, delineates a general abstract framework where
system theory ideas appear in disguised form. These system-theory
ingredients include: system matrix for an input/state/output linear
system, $Z$-transform  from a ``time domain'' to a ``frequency
domain'', and $Z$-transform of the output signal given by an observation
function applied to the initial condition plus a transfer function
applied to the $Z$-transform of the input signal.  Here we set down
the definitions and main results for the general Muhly-Solel
formalism and illustrate them for two specific types of multi-dimensional
linear systems:  (1) dissipative Fornasini-Marchesini state-space
representations with transfer function equal to a holomorphic
operator-valued function on the unit ball in ${\mathbb C}^d$, and
(2) noncommutative dissipative Fornasini-Marchsini  linear systems
with evolution along a free semigroup and with transfer function
defined on the noncommutative ball of   strict row contractions on a
Hilbert space.
\end{abstract}

\setcounter{equation}{0}
\section{Introduction}\label{intro}

In the last decade operator algebraists have worked on a
sophisticated  theory, based on the notions of $C^*$- and
$W^*$-correspondences, that encompasses many operator theoretical
results, for instance, from characteristic model and interpolation
theory \cite{MS05a,MS04}, both in time-varying and multivariable
settings, as well as for a range of more exotic examples. In the
recent paper \cite{MSSchur} a transfer function realization result
was obtained for the generalized Schur class (introduced in
\cite{MS05c}) in this $W^*$-correspondence setting; see also
\cite{BBFtH08}. It is the purpose of this note to enrich the
developed theory of \cite{MSSchur} with a system theory
interpretation, and show how the $W^*$-correspondences can be
specified to include some examples from multi-dimensional system theory.

In Section \ref{S:FM} we delineate the two examples that we
will obtain as special cases of the general Muhly-Solel framework:
(1) commutative Fornasini-Marchesini systems with evolution along the
positive-orthant integer lattice $\BZ_+^d$, and (2) noncommutative
Fornasini-Marchesini systems with evolution along the free semigroup
generated by $d$ letters. In Sections \ref{S:W*prelims} and \ref{S:CFpairs}
we lay out the somewhat formidable array of definitions needed to
formulate the general Muhly-Solel setup. The definition and basic properties
of the abstract input/state/output linear systems arising in the general
Muhly-Solel framework are given in Section \ref{S:systems}. Section
\ref{S:Point} is dedicated to the abstract point-evaluation and the
associated abstract $Z$-transform, while in Section \ref{S:Frequency
domain equations} we arrive at the frequency-domain version of the
system and basic properties of the associated observation function
and transfer function. An other type of linear systems which fits
the Muhly-Solel framework (1-D time-varying systems) and some other
types which do not are mentioned in the final Section \ref{S:other}.
Throughout the exposition, new ideas and concepts are illustrated
for the two basic examples as they are introduced.

The notation is mostly standard, but we mention here a few conventions.
For any index set $\Om$ we write $\ell^2(\Om)$ for Hilbert space of
complex-valued square summable functions over $\Om$:
\[
\left\{ f:\Om\to\BC \colon \mbox{$\sum_{\om\in\Om}|f(\om)|^2<\infty$} \right\},
\]
while $\ell^\infty(\Om)$ denotes the space of bounded complex-valued functions
over $\Om$. For a Hilbert space $\cH$ we use $\ell^2_\cH(\Om)$ and
$\ell^\infty_\cH(\Om)$ as short-hands for $\ell^2(\Om)\otimes\cH$ and
$\ell^\infty(\Om)\otimes\cH$, respectively.

\section{Two examples of multidimensional linear systems and their transfer functions} \label{S:FM}

In this section we review dissipative systems of Fornasini-Marchesini type; both commutative and
noncommutative. In a somewhat more general form the noncommutative case was studied in
\cite{BGM05,BGM06} with connections to robust control discussed in \cite{BGM06b}. See also
\cite{B07} for an overview of these types of systems and their applications.

\subsection{Commutative Fornasini-Marchesini systems}
Consider a system evolving over $\BZ_+^d$, the set of $d$-tuples of
nonnegative integers, given by the system equations:
\begin{equation}\label{CFMSysEqn}
\left\{\begin{array}{rcl}
x(n)\!\!\!\!&=&\!\!\!\! A_1 x(n-e_1)+\cdots +A_d x(n-e_d)\\
&&+B_1 u(n-e_1)+\cdots +B_d u(n-e_d),\\[.2cm]
y(n)\!\!\!\!&=&\!\!\!\! C x(n)+D u(n),
\end{array}\right.
\end{equation}
where $e_1,\ldots,e_d$ are the standard basis vectors in $\BR^d$ and
we set $x(n)=0$ for all $n\in\BZ^d$ that are not in $\BZ_+^d$, and
with a system matrix of the form
\begin{equation}\label{FMSysMat}
\mat{cc}{A_1&B_1\\\vdots&\vdots\\A_d&B_d\\C&D}:\mat{c}{\cX\\\cU}\to\mat{c}{\cX^d\\\cY}
\end{equation}
Here $A_j$, $B_j$, $C$ and $D$ are bounded linear operators and $\cX$, $\cU$ and $\cY$
are Hilbert spaces. We further assume the system matrix \eqref{FMSysMat} to be contractive.
In this case we speak of a dissipative system of commutative Fornasini-Marchasini type
(see \cite{FM}). Note that the output sequence $\{y(n)\}_{n\in\BZ_+^d}$ and state
sequence $\{x(n)\}_{n\in\BZ_+^d}$ are completely determined by the input sequences
$\{u(n)\}_{n\in\BZ_+^d}$ and the initial state $x(0)$. Moreover, the input, output and state
sequences satisfy the energy balance inequalities: For each $N\in\BZ_+$ we have
\[
\begin{array}{l}
\displaystyle
\sum_{n\in\BZ_+^d;|n|\leq N}\!\!\!\!\!\!\!\!\om(n)\|y(n)\|^2
+\!\!\!\!\!\!\!\!\sum_{n\in\BZ_+^d;|n|= N}\!\!\!\!\!\!\!\!\om(n)\|x(n)\|^2\leq\\[.2cm]
\displaystyle
\qquad\qquad\qquad\qquad\leq
\!\!\!\!\!\!\!\!\sum_{n\in\BZ_+^d;|n|\leq N}\!\!\!\!\!\!\!\!\om(n)\|u(n)\|^2 + \|x(0)\|^2
\end{array}
\]
where $\om(n)=\frac{|n|!}{n!}$ with $|n|=n_1+\cdots+n_d$ and $n!=n_1!\cdots n_d!$
for each $n=(n_1,\ldots,n_d)\in\BZ_+^d$. For a Hilbert space $\cV$ we write
$\ell^2_{\cV,\om}(\BZ_+^d)$ for the weighted $\ell^2$-space of sequences
$v=\{v(n)\}_{n\in\BZ_+^d}$ with $v(n)\in\cU$ and
\[
\|v\|^2=\sum_{n\in\BZ_+^d} \om(n)\|v(n)\|^2<\infty.
\]
It then follows from the energy balance inequalities that if the input sequence
$u$ is in $\ell^2_{\cU,\om}(\BZ_+^d)$, the output sequence $y$ is in
$\ell^2_{\cY,\om}(\BZ_+^d)$ and $\|y\|^2\leq\|u\|^2+\|x(0)\|^2$,
while the state sequence $x$ is bounded; more precisely
$\sup_{N\in\BZ_+}\sum_{|n|= N}\om(n)\|x(n)\|^2<\infty$. The state sequence
is in $\ell^2_{\cX,\om}(\BZ_+^d)$ if an additional stability assumption on
state operators $A_1, \dots, A_d$ is satisfied.

We write $\BB^d$ for the open unit ball in $\BC^d$:
\begin{equation}\label{BB}
\BB^d:=\left\{z=(z_1,\ldots,z_d)\in\BC^d \colon \sum_{j=1}^d|z_j|^2<1\right\}
\end{equation}
and for a $z=(z_1,\ldots,z_d)$ in $\BC^d$ we use the notation $z^n$
for $z_1^{n_1}\cdots z_d^{n_d}$ in case $n=(n_1,\ldots,n_d)\in\BZ_+^d$.
Given a bounded sequence $v=\{v(n)\}_{n\in\BZ_+^d}$, with $v(n)$ in a Hilbert
space $\cV$, we define a $\cV$-valued function $\widehat{v}$ on $\BB^d$
via the $Z$-transform:
\begin{equation}  \label{C-Ztrans}
\widehat{v}(z)=\sum_{n\in\BZ_+^d}z^nv(n).
\end{equation}
Provided that the input sequence is in $\ell^2_{\cU,\om}(\BZ_+^d)$,
the $Z$-transforms of the input, output and state sequences
are related through the frequency domain equations:
\begin{equation*}  
\left\{\begin{array}{rcl}
\widehat{x}(z)\!\!\!\!&=&\!\!\!\! x(0)+\!\!\left(\sum_{k=1}^d z_kA_k\right)\!\widehat{x}(z)+\!\!
\left(\sum_{k=1}^d z_kB_k\right)\!\widehat{u}(z),\\[.2cm]
\widehat{y}(z)\!\!\!\!&=&\!\!\!\! C\widehat{x}(z)+D\widehat{u}(z),
\end{array}\right.
\end{equation*}
from which one finds that the input/inital state to output map is given by
\begin{equation}\label{Cinput-output}
\widehat{y}(z)=F(z)\widehat{u}(z)+W(z)x(0),
\end{equation}
where $F$ and $W$ are the transfer function and observation function given by
\begin{eqnarray}
\notag
F(z)&=&D+C(I-z_1A_1-\cdots -z_dA_d)^{-1}\times\\
 \label{CFMtrans-obs}
&&\qquad\qquad\times(z_1B_1+\cdots+z_dB_d),\\
\notag
W(z)&=&C(I-z_1A_1-\cdots -z_dA_d)^{-1}.
\end{eqnarray}

\subsection{Noncommutative Fornasini-Marchesini systems}  \label{S:IntroNCFM}
We write $\cF_d$ for the free semigroup generated by the $d$-letter
alphabet $\{1,\ldots,d\}$; elements of $\cF_d$ are words of the form
$\al=i_{N}\cdots i_2 i_1$ with $i_j\in\{1,\ldots,d\}$ for
$j=1,\ldots,N$. We also include the empty word $\emptyset$ that
serves as the unit element with respect to the semigroup
action: concatenation of words. With $|\al|$ we indicate the
length of the word $\al$.

Now consider a system that evolves along the free semigroup $\cF_d$ given by the equations
\begin{equation}\label{NCFMSysEqn}
\left\{\begin{array}{rcl}
x(k\cdot\al)\!\!\!\!&=&\!\!\!\! A_kx(\al) +B_k u(\al),\ k=1,\ldots,d,\\[.2cm]
y(\al)\!\!\!\!&=&\!\!\!\! C x(\al)+D u(\al),
\end{array}\right.
\end{equation}
where we have a system matrix of the form \eqref{FMSysMat} that we
assume to be contractive. For notational simplicity we  assume
$\cU=\cY=\BC$. Again, the output sequence $\{y(\al)\}_{\al\in\cF_d}$
and state sequence $\{x(\al)\}_{\al\in\cF_d}$ are completely
determined by the input sequence $\{u(\al)\}_{\al\in\cF_d}$ and an
initial condition: $x(\emptyset)=x_\emptyset$, and in this case we
have energy balance inequalities of the form
\[
\begin{array}{l}
\displaystyle
\sum_{\al\in\cF_d;|\al|\leq N}\!\!\!\!\!\!\!\!\|y(\al)\|^2
+\!\!\!\!\!\!\!\!\sum_{\al\in\cF_d;|\al|= N}\!\!\!\!\!\!\!\!\|x(\al)\|^2\leq\\[.2cm]
\displaystyle
\qquad\qquad\qquad\qquad\leq
\!\!\!\!\!\!\!\!\sum_{\al\in\cF_d;|\al|\leq N}\!\!\!\!\!\!\!\!\|u(\al)\|^2 + \|x(\emptyset)\|^2
\end{array}
\]
for each $N\in\BZ_+$. Hence an input sequence $u$ from $\ell^2(\cF_d)$ produces an output
sequence $y$ in $\ell^2(\cF_d)$ so that $\|y\|^2\leq\|u\|^2+\|x(\emptyset)\|^2$, and the
state sequence $x$ is bounded; more precisely $\sup_{N\in\BZ_+}\sum_{|\al|= N}\|x(\al)\|^2<\infty$.
In case the state operators $A_1, \dots, A_d$ satisfy an additional stability the
state sequence $x$ is in $\ell^2_\cX(\cF_d)$.

For this example we shall use a frequency domain consisting of
operator-tuples rather than  tuples of complex numbers.  We
therefore  assume that we are given an additional Hilbert space
$\cK$  and a  $d$-tuple $T=(T_1,\ldots,T_d)$ of operators on $\cK$.
Such a $d$-tuple $T$  is said to be a  {\em strict row contraction}
if the operator matrix $\mat{ccc}{T_1&\cdots&T_d}$ defines a strict
contraction from $\cK^d$ into $\cK$. For a strict row contraction
$T$ we define the noncommutative functional calculus
\begin{equation}\label{NonComFunCalc}
T^\al=T_{i_N}\cdots T_{i_1}\text{ if }\al=i_{N}\cdots i_1,\quad T^\emptyset=I_\cK.
\end{equation}
Then for $\cW=\BC,\cX$ and $w=\{w(\al)\}_{\al\in\cF_d}$ a bounded $\cW$-valued
sequence we define a function $\widehat{w}$ on the set of strict row contractions
 via the noncommutative $Z$-transform:~For each strict row contraction $T$ we set
\begin{equation}   \label{NC-Ztrans}
\widehat{w}(T)=\sum_{\al\in\cF_d}w(\al)\otimes T^\al\in\cL(\cK,\cW\otimes\cK).
\end{equation}
Then the system \eqref{NCFMSysEqn} is
described in the frequency domain by the equations
\begin{equation}\label{NCFMSysEqnFD}
\left\{\begin{array}{rcl}
\widehat{x}(T)\!\!\!\!&=&\!\!\!\! x(\emptyset)+\left(\sum_{k=1}^d A_k\otimes T_k\right)\widehat{x}(T)+\\
&&\qquad\qquad+\left(\sum_{k=1}^d B_k\otimes T_k\right)\widehat{u}(T),\\[.2cm]
\widehat{y}(T)\!\!\!\!&=&\!\!\!\! (C\otimes I_\cK)\widehat{x}(T)+(D\otimes I_\cK)\widehat{u}(T),
\end{array}\right.
\end{equation}
and the input/intial state-to-output map is given by
\begin{equation}\label{ii-o}
\widehat{y}(T)=F(T)\widehat{u}(T)+W(T)x(\emptyset),
\end{equation}
where we define the transfer function $F$ and observation function $W$ on the
set of strict row contractions $T$ by
\begin{eqnarray}
\notag
F(T)&=&D+C(I-A_1\otimes T_1-\cdots -A_d\otimes T_d)^{-1}\times\\
\notag
&&\times (B_1\otimes T_1+\cdots+B_d\otimes T_d)\\
 \label{NCtrans-obs}
&=&D+\sum_{k=1}^d\sum_{\al\in\cF_d}CA^\al B_k T^{\al\cdot k}\\
\notag
W(T)&=&C(I-A_1\otimes T_1-\cdots -A_d\otimes T_d)^{-1}\\
\notag
&=&\sum_{\al\in\cF_d}(CA^\al)\otimes T^\al.
\end{eqnarray}

It is possible to symmetrize the noncommutative Fornasini-Marchesini
system equations \eqref{NCFMSysEqn} to recover the commutative
Fornasini-Marchesini system equations \eqref{CFMSysEqn}.  To
describe this in detail, we introduce the abelianization map
\begin{equation}\label{abel}
\Abel(\al)=(n_1,\ldots,n_d)\quad (\al\in\cF_d)
\end{equation}
where, for $j=1,\ldots,d$, the integer $n_j$ is equal to the number
of times the letter $j$ appears in the word $\al$.
Let $\{u(\alpha), x(\alpha), y(\alpha)  \colon \alpha \in
\cF_d    \}$ be a system
trajectory satisfying \eqref{NCFMSysEqn}, for each $n \in {\mathbb
Z}^d_+$ define $\overline{u}(n) = \sum_{\alpha \colon \Abel(\alpha)
= n} u(\alpha)$, $\overline{x}(n) = \sum_{\alpha \colon
\Abel(\alpha) = n} x(\alpha)$, $\overline{y}(n) = \sum_{\alpha
\colon \Abel(\alpha) = n} y(\alpha)$.
Then one can check that the
more detailed system equations \eqref{NCFMSysEqn} imply that the
aggregate quantities $\{ \overline{u}(n), \overline{x}(n),
\overline{y}(n) \colon n \in {\mathbb Z}^d_+\}$ satisfy the system
equations \eqref{CFMSysEqn}.  Moreover, in general if the $d$-tuple
$T = (T_1, \dots, T_d)$ of operators on $\cK$ is commutative, then
$T^\alpha = T^{\alpha^\prime}=: T^n$ whenever
$\Abel(\alpha)=\Abel(\alpha')=n \in {\mathbb Z}^d_+$ and
the value of the $Z$-transform $\widehat w(T)$ in \eqref{NC-Ztrans}
can be expressed as $\widehat w(T) = \sum_{n \in   {\mathbb Z}^d_+}
\overline{w}(n) T^n$.  This applies in particular to the case where
we take ${\mathcal K} = {\mathbb C}$ to be one-dimensional; then $T
= (T_1, \dots, T_d)$ becomes a $d$-tuple of scalars $z = (z_1,
\dots, z_d) \in {\mathbb B}^d$, and the noncommutative $Z$-transform
\eqref{NC-Ztrans} collapses to the commutative $Z$-transform
\eqref{C-Ztrans} (with $\overline{w}$ in place of $v$).

\section{$C^*$-correspondence preliminaries} \label{S:W*prelims}

Let $\cA$ and $\cB$ be $C^*$-algebras and $E$ a linear space. We say
that $E$ is an {\em $(\cA,\cB)$-correspondence} when $E$ is a
bi-module with respect to a right $\cB$-action and a left
$\cA$-action, and $E$ is endowed with a $\cB$-valued inner product
$\inn{\ }{\ }_E$ satisfying the following axioms: For any
$\la,\mu\in\C$, $\xi,\eta,\ze\in E$, $a\in\cA$ and $b\in\cB$, it is the case that
\begin{enumerate}
\item $\inn{\la\xi+\mu\ze}{\eta}_E =\la\inn{\xi}{\eta}_E+\mu\inn{\ze}{\eta}_E;$
\item $\inn{\xi\cdot b}{\eta}_E=\inn{\xi}{\eta}_E b;$
\item $\inn{a\cdot\xi}{\eta}_E=\inn{\xi}{a^*\cdot\eta}_E;$
\item $\inn{\xi}{\eta}_E^*=\inn{\eta}{\xi}_E;$
\item $\inn{\xi}{\xi}_E\geq 0 \ (\text{in }\cB);$
\item $\inn{\xi}{\xi}_E=0\mbox{ implies that }\xi=0$,
\end{enumerate}
and that $E$ is a Banach space with respect to the norm
\[
\|\xi\|_E=\|\inn{\xi}{\xi}_E\|^\half_\cB\quad(\xi\in E),
\]
where $\|\ \|_\cB$ denotes the norm of $\cB$. We also impose that
\begin{equation*}
(\la\xi)\cdot b=\xi\cdot(\la b)\ands (\la
a)\cdot\xi=a\cdot(\la\xi)
\end{equation*}
for any $\la\in\C,a\in\cA,b\in\cB,\xi\in E$.
In practice we usually write $\inn{\ }{\ }$ and $\|\ \|$ for the
inner product and norm on $E$, and in case $\cA=\cB$ we say that $E$
is an $\cA$-correspondence.

\begin{example}\label{E:1}{\em
{\bf Part A:}  For  the commutative Fornasini-Marchensini system interpretation,
we take $\cA=\cL(\cU)$ and $E=\text{col}_{j=1}^d[\cL(\cU)]=\cL(\cU,\cU^d)$, where $\cU$ is
a given Hilbert space. Then $E$ can be seen as an $\cA$-correspondence with left and right
$\cA$-action and inner product given by
\begin{align*}
a \cdot \text{col}_{j=1}^d[T_j] \cdot a'=\text{col}_{j=1}^d[aT_ja'],\\
\inn{\text{col}_{j=1}^d[T_j]}{\text{col}_{j=1}^d[R_j]}=\sum_{j=1}^d R_j^*T_j.
\end{align*}
{\bf Part B:} For the noncommutative Fornasini-Marchensini system interpretation we take
the same correspondence, but then with $\cU=\BC$.
}\end{example}

Given two $(\cA,\cB)$-correspondences $E$ and $F$, the set of
bounded linear operators from $E$ to $F$ is denoted by $\cL(E,F)$.
It may happen that a $T\in\cL(E,F)$ is not adjointable, i.e., it is
not necessarily the case that
\[
\inn{T\xi}{\ga}=\inn{\xi}{T^*\ga}\quad(\xi\in E,\ga\in F)
\]
for some $T^*\in\cL(F,E)$.  We use $\cL^a(E,F)$ to indicate the set of operators
in $\cL(E,F)$ which are adjointable. As usual we have the
abbreviations $\cL(E)$ and $\cL^a(E)$ in case $F=E$.

The third inner-product axiom implies that the left $\cA$-action can
be identified with a $*$-representation (meaning a nondegenerate
$*$-homomorphism) $\vph$ of $\cA$ into the
$C^*$-algebra $\cL^a(E)$. In case this $*$-representation $\vph$ is
specified we will occasionally write $\vph(a) \xi$ instead of
$a\cdot \xi$.

Furthermore, an operator $T\in\cL(E,F)$ is said to be a {\em right
$\cB$-module map} if
\[
T(\xi\cdot b)=T(\xi)\cdot b\quad(\xi\in E, b\in\cB),
\]
and a {\em left $\cA$-module map} whenever
\[
T(a\cdot\xi)=a\cdot T(\xi)\quad(\xi\in E, a\in\cA).
\]

We will have a need for various constructions which create new
correspondences out of given correspondences.

Given two $(\cA, \cB)$-correspondences $E$ and $F$, we define the
{\em direct-sum correspondence} $E \oplus F$ to be the direct-sum
vector space $E \oplus F$ together with the left
$\cA$-action and right $\cB$-action and the direct-sum $\cB$-valued
inner-product defined by setting for each $\xi,\xi'\in E$,
$\ga,\ga'\in F$, $a\in\cA$ and $b\in\cB$:
\[
\begin{array}{c}
a\cdot(\xi\oplus\ga)=(a\cdot\xi) \oplus (a\cdot\ga),\\[.1cm]
(\xi\oplus\ga)\cdot b=(\xi\cdot b)\oplus(\ga\cdot b),\\[.1cm]
\inn{\xi\oplus\ga}{\xi'\oplus\ga'}_{E\oplus F}=
\inn{\xi}{\xi'}_E+\inn{\ga}{\ga'}_F.
\end{array}
\]
Bounded linear operators between direct-sum correspondences admit
operator matrix decompositions in precisely the same way as in the
Hilbert space case ($\cB=\C$), while adjointability and the left and
right module map property of such an operator corresponds to the
operators in the decomposition being adjointable, or left or right
module maps, respectively.

Now suppose we are given three $C^{*}$-algebras $\cA, \cB$ and
${\mathcal C}$, an $(\cA, \cB)$-correspondence $E$ and
a $(\cB,\cC)$-correspondence $F$. Then we define the {\em
tensor-product correspondence} $E \otimes F$ to be the completion of
the linear span of all tensors $\xi\otimes\ga$ (with $\xi\in E$ and
$\ga\in F$) subject to the identification
\begin{equation}\label{balance}
 (\xi\cdot b)\otimes\ga=\xi\otimes(b\cdot\ga)\quad(\xi\in E,\ga\in F,b\in\cB),
\end{equation}
with left $\cA$-action, right $\cC$-action and the $\cC$-valued
inner-product defined, initially just for pure tensors, by setting for
each $\xi,\xi'\in E$, $\ga,\ga'\in F$, $a\in\cA$ and $c\in\cC$:
\[
\begin{array}{c}
a\cdot(\xi\otimes\ga)=(a\cdot\xi)\otimes\ga,\qquad
(\xi\otimes\ga)\cdot c=\xi\otimes(\ga\cdot c),\\[.1cm]
\inn{\xi\otimes\ga}{\xi'\otimes\ga'}_{E\otimes F}=
\inn{\inn{\xi}{\xi'}_E\cdot\ga}{\ga'}_F,
\end{array}
\]
and then extending by linearity and continuity to $E\otimes F$.
In case the left action on $F$ is given by the $*$-representation
$\vph$ we occasionally emphasize this by writing $E\otimes_\vph F$
for $E\otimes F$.
It is more complicated to characterize the bounded linear operators
between tensor-product correspondences. One way to construct such
operators is as follows.  Let $E$ and $E'$ be
$(\cA,\cB)$-correspondences and $F$ and $F'$
$(\cB,\cC)$-correspondences. Furthermore, let $X\in\cL(E,E')$ be a right module map and
$Y\in\cL(F,F')$ a left module map. We then write $X\otimes Y$ for
the {\em tensor product operator} in $\cL(E\otimes F,E'\otimes F')$
determined by
\begin{equation}\label{tensoroperator}
X\otimes Y (\xi\otimes\ga)=(X\xi)\otimes(Y\ga)\quad
(\xi\otimes\ga\in E \otimes F).
\end{equation}
If, in addition, $X$ is also a left module map, then $X\otimes Y$ is
a left module map, while $Y$ also being a right module map
guarantees that $X\otimes Y$ is a right module map. Moreover,
$X\otimes Y$ is adjointable when $X$ and $Y$ are both adjointable,
and in that case $(X\otimes Y)^*=X^*\otimes Y^*$.
Notice that the left action on $E\otimes F$ can now be written as
$a\mapsto\varphi(a)\otimes I_F\in\cL^a(E\otimes F)$, where
$I_F$ denotes the identity operator on $F$ and $\vph$ indicates
the left $\cA$-action on $E$.

\section{Correspondence-representation pairs and their duals}\label{S:CFpairs}

In the remainder of the paper we consider the situation where $\cA = \cB$,
i.e, $E$ is an $\cA$-correspondence. We also restrict our attention
to the case where $\cA$ is a von Neumann algebra and $E$ an
$\cA$-$W^{*}$-correspondence. This means that $E$ is an
$\cA$-correspondence which is also {\em self-dual} in the sense that
any right $\cA$-module map $\rho \colon E \to \cA$ is given by
taking the inner-product against some element $e_{\rho}$ of $E$:
\begin{equation}  \label{rho}
\rho(e) = \langle e, e_{\rho} \rangle_{E} \in \cA.
\end{equation}
Suppose that we are also given an
auxiliary Hilbert space $\cE$ and a $*$-representation $\sigma
\colon \cA \to \cL(\cE)$; as this will be the setting for much of
the analysis to follow, we refer to such a pair $(E, \sigma)$ as a
{\em Correspondence-Representation pair}, or {\em CR-pair} for
short. We further assume that $\si$ is faithful (injective) and
normal ($\si$-weakly continuous). The Hilbert space $\cE$
becomes an $(\cA,\C)$-correspondence with left $\cA$-action given by
$\sigma$:
$$
 a \cdot y = \sigma(a) y \text{ for all } a \in \cA \text{ and } y
 \in \cE.
 $$
Thus we can form the tensor-product $(\cA, {\mathbb
C})$-correspondence $E \otimes_{\sigma} \cE$. We then write $E^{\sigma}$ for
the set of bounded linear operators $\mu\colon \cE\to E
\otimes_{\sigma} \cE$ which are also left $\cA$-module maps:
\begin{equation} \label{Esigma0}
E^{\sigma} = \{ \mu\colon \cE\to E \otimes_{\sigma} \cE \colon \mu
\sigma(a)  = ( \varphi(a) \otimes I_{\cE}) \mu \}.
\end{equation}
It turns out that $E^{\sigma}$ is itself a $W^*$-correspondence, not
over $\cA$ but over the $W^*$-algebra
$$
\sigma(\cA)' = \{b \in \cL(\cE) \colon b \sigma(a) = \sigma(a) b
\text{ for all } a \in \cA\} \subset \cL(\cE)
$$
(the {\em commutant} of the image $\sigma(\cA)$ of the
$*$-representation $\sigma$ in $\cL(\cE)$) with left and right
$\sigma(\cA)'$-action and $\sigma(\cA)'$-valued inner-product
$\langle \cdot, \cdot \rangle_{E^{\sigma}}$ given by
\[
\mu\cdot b=\mu b,\quad b\cdot\mu=(I_E\otimes b)\mu,\quad
\langle \mu, \nu \rangle_{E^{\sigma}} = \nu^{*}\mu.
\]
Let $\iota:\si(\cA)'\to\cL(\cE)$ be
the identity $*$-representation: $\iota(b)=b$ for each
$b\in\si(\cA)'$. Then $(E^\si,\iota)$ is a CR-pair that is said to
be {\em dual} to $(E,\si)$; see Section 3 in \cite{MS04} for further
details. In particular, $\cE$ can be seen both as an
$(\cA,\BC)$-correspondence (with left action given by $\si$) and a
$(\si(\cA)',\BC)$-correspondence (with left action given by
$\iota$).

\begin{example}[Continuation of Example \ref{E:1}]\label{E:2}{\em $ $\\
{\bf Part A:}
We take $\cE=\cU$ and for the $*$-representation $\si$ the
identity map: $\si(T)=T$ for each $T\in\cL(\cU)$. Then
$E\otimes\cE=\cU^d$, $\si(\cA)'=\BC I_\cU:=\{\mu I_\cU\colon \mu\in\BC\}$ and
\begin{align*}
E^\si=\{\col_{j=1}^d[\mu_j I_\cU]\colon \mu_j\in\BC\}\subset\cL(\cU,\cU^d).
\end{align*}
The spaces $\si(\cA)'$ and $E^\si$ shall be identified with $\BC$ and $\BC^d$, respectively.\\
{\bf Part B:} Let $\cK$ be some Hilbert space. We then take $\cE=\cK$ and
for $\si:\BC\to\cL(\cK)$ the $*$-representation $\si(\la)=\la I_\cK$ for
each $\la\in\BC$. Then $E\otimes\cE=\cK^d$, $\si(\cA)'=\cL(\cK)$
and
\[
E^\si=\{\col_{j=1}^d[T_j]\colon T_j\in\cL(\cK)\}=\cL(\cK,\cK^d).
\]

In these two examples the duality mentioned above takes the following form.
In Part A we start with $E=\cL(\cU,\cU^d)$, $\cA=\cL(\cU)$ and $\si$ the identity map
on $\cL(\cU)$. We then find the dual CR-pair $(E^\si,\iota)$ with $E^\si=\BC^d$,
$\si(\cA)'=\BC$ (after the identification) and $\iota(\la)=\la I_\cU$. This is precisely
the CR-pair we use in Part 2 (if we take $\cK=\cU$). It then follows from Part 2 that
forming the CR-pair dual to $(E^\si,\iota)$ we return to $(E,\si)$. In this sense the
two examples are also each others dual.
}\end{example}

Since $E$ is an $\cA$-correspondence, we can inductively define the
self-tensor-product $\cA$-correspondences $E^{\otimes n}$ via
$E^{\otimes 1}=1$ and $E^{\otimes n+1} = E \otimes (E^{\otimes n})$
for $n = 1,2,\dots$. Formally we set $E^{\otimes 0} = \cA$. Using
the notation in \eqref{tensoroperator}, the left $\cA$-action on
$E^{\otimes n}$ can be written as $\varphi(a)\otimes I_{E^{\otimes
n-1}}$ if we use $\varphi:\cA\to\cL^{a}(E)$ to indicate the left
$\cA$-action of $E$. We may then also form the $(\cA,\BC)$-correspondences
$E^{\otimes n}\otimes_\si\cE$. Besides the left $\cA$-action, these Hilbert
spaces also admit a left $\si(\cA)'$-action via the $*$-representation
$b\mapsto I_{E^{\otimes n}}\otimes b$ for $b\in\si(\cA)'$. Furthermore, for
a $\mu\in E^\si$ we can form the tensor product operator
$I_{E^{\otimes n}}\otimes\mu$ that maps $E^{\otimes n}\otimes\cE$ into
$E^{\otimes n+1}\otimes\cE$. In a similar way we define the
$\si(\cA)'$-correspondences $(E^{\si})^{\otimes n}$ and
$(\si(\cA)',\BC)$-correspondences $(E^{\si})^{\otimes n}\otimes_\iota\cE$
for $n\in\BZ_+$.

One of the consequences of the duality between $(E,\si)$ and
$(E^\si,\iota)$ is the following result from \cite{MS04}.

\begin{proposition}\label{P:duality}{\em
For each $n\in\BZ_+$ the map
$\Phi_n:(E^\si)^{\otimes n}\otimes_\iota\cE\to E^{\otimes n}\otimes_\si\cE$ induced
by the identity
\[
\Phi_n(\mu_n\otimes\cdots\otimes\mu_1\otimes e)
=(I_{E^{\otimes n-1}}\otimes\mu_n)\cdots(I_E\otimes\mu_2)\mu_1e
\]
is a unitary left $\si(\cA)'$-module map.}
\end{proposition}

\begin{example}[Continuation of Example \ref{E:2}]\label{E:3}{\em $ $\\
{\bf Part A:} For $n\in\BZ_+$ set $d_n=\#\{\al\in\cF_d,|\al|=n\}=d^n$. Then
\[
E^{\otimes n}=\cL(\cU,\cU^{d_n})\ands E^{\otimes n}\otimes_\si\cE=\cU^{d_n},
\]
while
\[
(E^\si)^{\otimes n}=\cL(\BC,\BC^{d_n})=\BC^{d_n}\ands
(E^\si)^{\otimes n}\otimes_\iota\cE=\cU^{d_n}.
\]
In particular, $E^{\otimes n}\otimes_\si\cE$ and $(E^\si)^{\otimes n}\otimes_\iota\cE$
are the same Hilbert space, and moreover, the left $\si(\cA)'=\BC$-actions on both
spaces is just scalar multiplication. The left $\cA=\cL(\cU)$-action on $\cU^{d_n}$
is given by
\begin{equation}\label{ExLeftAction}
T\left( \oplus_{j=1}^{d_n} e_j \right)= \oplus_{j=1}^{d_n} Te_j,
\end{equation}
for $T\in\cL(\cU)$ and $e_j\in\cU$. Moreover, the map $\Phi_n$
defined in Proposition \ref{P:duality} is the identity map on $\cU^{d_n}$.\\
{\bf Part B:} Due to the duality between the two examples, as observed in
Example \ref{E:2}, we find:
\begin{align*}
E^{\otimes n}=\cL(\BC,\BC^{d_n})=\BC^{d_n}\ands
E^{\otimes n}\otimes_\si\cE=\cK^{d_n},\\
(E^\si)^{\otimes n}=\cL(\cK,\cK^{d_n})\ands
(E^\si)^{\otimes n}\otimes_\iota\cE=\cK^{d_n},
\end{align*}
with left $\si(\cA)'=\cL(\cK)$-action on $\cK^{d_n}$ given by \eqref{ExLeftAction},
with $T\in\cL(\cK)$, $e_j\in\cK$. Moreover, $\Phi_n$ is again the identity map on $\cK^{d_n}$.
}\end{example}

\section{Linear dissipative discrete-time systems}\label{S:systems}

Let $(E,\si)$ be a CR-pair as in Section \ref{S:CFpairs}. A
{\em (linear) system} associated with the pair
$(E,\si)$ is a set $\Sys$ consisting of four left $\si(\cA)'$-module maps:
\begin{equation}\label{SysMat}
\mat{cc}{A&B\\C&D}:\mat{c}{\cX\\\cE}\to\mat{c}{E^\si\otimes\cX\\\cE},
\end{equation}
and a $*$-representation $\pi:\si(\cA)'\to\cL(\cX)$ that provides a
left $\si(\cA)'$-action on the Hilbert space $\cX$, called the {\em state space},
making $\cX$ into a $(\si(\cA)',\BC)$-correspondence.
In this paper we are primarily interested in the case that the {\em system
matrix} \eqref{SysMat} is contractive, in which case we speak of a
{\em dissipative system}. For simplicity we only consider here the case that
the {\em input space} and {\em output space} are the same, i.e., both equal
to $\cE$; various techniques exist to extend to the case with different input
and output spaces, cf., \cite{MSSchur,BBFtH08}.

Given a system $\Sys$, for each $n\in\BZ_+$ we define the operators:
\begin{equation*}
\begin{array}{c}
A_n=I_{(E^\si)^{\otimes n}}\otimes A,\quad
B_n=(I_{(E^\si)^{\otimes n}}\otimes B)\Phi_n^*\\[.2cm]
C_n=\Phi_n(I_{(E^\si)^{\otimes n}}\otimes C),\quad
D_n=\Phi_n(I_{(E^\si)^{\otimes n}}\otimes D)\Phi_n^*,
\end{array}
\end{equation*}
where $\Phi_n$ is as defined in Proposition \ref{P:duality},
and associate with $\Sys$ the following linear system of equations:
\begin{equation}\label{SysEqn}
\Si=\left\{\begin{array}{rcl}
x(n+1)\!\!\!\!&=&\!\!\!\! A_n x(n)+B_nu(n),\\[.2cm]
y(n)\!\!\!\!&=&\!\!\!\! C_n x(n)+D_n u(n),
\end{array}\right.(n\in\BZ_+)
\end{equation}
where $u(n),y(n)\in E^{\otimes n}\otimes\cE$, and $x(n)\in (E^\si)^{\otimes n}\otimes\cX$.
Note that an {\em input sequence} $\{u(n)\}_{n\in\BZ_+}$ and initial state $x(0)$ uniquely
determine the {\em output sequence} $\{y(n)\}_{n\in\BZ_+}$ and {\em state sequence}
$\{x(n)\}_{n\in\BZ_+}$.

\begin{example}[Continuation of Example \ref{E:3}]\label{E:4}{\em $ $\\
{\bf Part A:} Let $\Sys$ be a dissipative system associated with the CR-pair
of Example \ref{E:2}, Part A. Then $\pi:\BC\to\cL(\cX)$ must be equal to
$\pi(\la)=\la I_\cX$, and thus $E^\pi\otimes\cX=\cX^d$. The system then has the form
\begin{equation}\label{ExSysmat}
\mat{cc}{A&B\\C&D}=\mat{c|c}{A_1&B_1\\\vdots&\vdots\\A_d&B_d\\\hline C&D}:
\mat{c}{\cX\\\cU}\to\mat{c}{\cX^d\\\cU},
\end{equation}
and, since $\si(\cA)'=\BC$, the $\si(\cA)'$-module map properties are satisfied trivially.
The input and output sequences ${\bf u}=\{\widetilde{u}(n)\}_{n\in\BZ_+}$ and
${\bf y}=\{\widetilde{y}(n)\}_{n\in\BZ_+}$ have elements $\widetilde{u}(n)$ and
$\widetilde{y}(n)$ in $\cU^{k_n}$, while the elements of the state sequence
${\bf x}=\{\widetilde{x}(n)\}_{n\in\BZ_+}$ are in $\cX^{k_n}$. The system equations
take the form
\begin{equation}\label{ExSysEq}
\Si:=\left\{\!\!
\begin{array}{rcl}
\widetilde{x}(n+1)\!\!\!\!&=&\!\!\!\!\tilde A_n \widetilde{x}(n)+\tilde B_n\widetilde{u}(n),\\[.2cm]
\widetilde{y}(n)\!\!\!\!&=&\!\!\!\!\tilde C_n \widetilde{x}(n)+\tilde D_n\widetilde{u}(n),
\end{array}\ (n\in\BZ_+)\right.
\end{equation}
with
\begin{equation*}
\begin{array}{l}
\tilde A_n=\operatorname{diag}_{i=1}^{d_n}\left(\mat{c}{A_1\\\vdots\\A_d}\right),\ \
\tilde B_n=\operatorname{diag}_{i=1}^{d_n}\left(\mat{c}{B_1\\\vdots\\B_d}\right),\\[.2cm]
\tilde C_n=\operatorname{diag}_{i=1}^{d_n}(C),\quad
\tilde D_n=\operatorname{diag}_{i=1}^{d_n}(D),
\end{array}
\end{equation*}
where $\operatorname{diag}$ forms a block diagonal operator, i.e.,
$\tilde A_n\in\cL(\cX^{d_n},\cX^{d_{n+1}}),$ $\tilde B_n\in\cL(\cU^{d_n},\cX^{d_{n+1}}),$ etc.
After making the identification of an element $\widetilde{u}(n)$  in $\cU^{d_n}$ with a set of
elements of $\cU$ of the form $\{u(\al)\colon \al\in\cF_d,\, |\al|=n\}$ via the bijective map
\begin{eqnarray}
\label{nu}
&\nu_n:\{\al\in\cF_d \colon |\al|=n \}\to \{1,\ldots,d_n\},&\\
\nonumber
&\nu_n(\al)=1+\sum_{j=1}^n(i_j-1)d^{j-1}\quad \text{if}\quad\al=i_n\cdots i_1,&
\end{eqnarray}
and similarly for $\cX^{d_n}$, and after unpacking the system equations \eqref{ExSysEq}, we
arrive at the noncommutative system
\begin{equation}\label{NCFMSysEqn2}
\left\{\begin{array}{rcl}
x(k\cdot\al)\!\!\!\!&=&\!\!\!\! A_k x(\al) +B_k u(\al),\ k=1,\ldots,d,\\[.2cm]
y(\al)\!\!\!\!&=&\!\!\!\! C x(\al)+D u(\al),
\end{array}\right.
\end{equation}
that evolves over $\cF_d$.
We obtain the commutative system \eqref{CFMSysEqn}  after applying the abelianization procedure.\\
{\bf Part B:} In this case the system matrix of a dissipative system $\Sys$ associated
with the CR-pair of Example \ref{E:2}, Part B, is as in \eqref{ExSysmat}, but with $\cU=\cK$,
and $\pi$ is a $*$-representation from $\cL(\cK)$ into $\cL(\cX)$. Moreover, the operators in
$\Sys$ satisfy the following intertwining relations:
\begin{equation}\label{intertwine}
\begin{array}{c}
\pi(T)A_j=A_j\pi(T),\quad \pi(T)B_j=B_jT,\\[.2cm]
TC=C\pi(T),\quad TD=DT,
\end{array}
\end{equation}
for any $T\in\cL(\cK)$ and $j=1,\ldots,d$. In particular, $D$ is a
scalar multiple of $I_\cK$. Assume that $\pi$ is a scalar multiple of the
identity representation; i.e., there exists a Hilbert space $\bar{\cX}$ so that
$\cX=\bar{\cX}\otimes\cK$ and $\pi(T)=I_{\bar{\cX}}\otimes T$ for each $T\in\cL(\cK)$.
In fact, if $\cK$ is finite dimensional, then any representation $\pi$ has this form
(see \cite{Arveson76}).
Then \eqref{intertwine} simplifies to
\begin{equation}\label{inertwine2}
A_j=\bar{A}_j\otimes I_\cK,\quad B_j=\bar{B}_j\otimes I_\cK,\quad C=\bar{C}\otimes I_\cK,
\end{equation}
with $\bar{A}_j\in\cL(\bar{\cX})$, $\bar{B}_j\in\cL(\BC,\bar{\cX})$ and $\bar{C}\in\cL(\bar{\cX},\BC)$
for $j=1,\ldots,d$.
Due to this special form of the system operators, we see that, if the input signal
$\{u(\alpha)\}_{\alpha \in \cF_d}$, identified with the actual input signal
$\{\widetilde{u}(n)\}_{n\in\BZ_+}$ as in Part A, and initial state $x(0)$ are
restricted to have the form $\{\ov{u}(\al)k_0\}_{\al\in\cF_d}$ and $\ov{x}(0)\otimes k_0$,
respectively, where $\ov{u}(\al)$ is scalar-valued, $\ov{x}(0)\in\bar{\cX}$ and $k_0$ is a
fixed vector in $\cK$, then the output sequence $\{y(\alpha)\}_{\alpha \in \cF_d}
= \{ \overline{y}(\alpha) k_0\}_{\alpha \in \cF_d}$ will have the same form as the input
sequence, i.e., $\overline{y}(\alpha)\in\BC$, while the state sequence is of the form
$\{x(\alpha)\}_{\alpha \in \cF_d}= \{ \overline{x}(\alpha)\otimes k_0\}_{\alpha \in \cF_d}$
with $\overline{x}(\alpha)\in\overline{\cX}$.
{}From now on, we shall restrict the inputs and initial states to vectors of this
form for a fixed vector $k_0 \in \cK$. In that case, the transformation to the
noncommutative system equations occurs in the same way as in Part A,
and we arrive at \eqref{NCFMSysEqn} with $A_j,$ $B_j$ and $C$ replaced with
$\bar{A}_j$, $\bar{B}_j$ and $\bar{C}$.
}\end{example}

We now define $\cF^2(E)=\oplus_{n=0}^\infty E^{\otimes n}$, which is an $\cA$-correspondence
with left $\cA$-action given by the $*$-representation $\vph_\infty:\cA\to\cL^a(\cF^2(E))$
that is defined via
\[
\vph_\infty(a)=\diag_{n\in\BZ_+}(\vph(a)\otimes I_{E^{\otimes n-1}})
\]
where $\vph$ is used to indicate the $*$-representation that provides the left
$\cA$-action on $E$. Then we can form the tensor product space
$\cF^2(E,\si)=\cF^2(E)\otimes_\si\cE$, which in turn is an $(\cA,\BC)$-correspondence.
We refer to $\cF^2(E)$ and $\cF^2(E,\si)$ as the {\em Fock spaces} associated with the
correspondence $E$ and the CR-pair $(E,\si)$, respectively.

Besides the left $\cA$-action, the Hilbert space $\cF^2(E,\si)$ is also equipped with
a left $\si(\cA)'$-action that is given by the $*$-representation
$b\mapsto I_{\cF^2(E)}\otimes b\in\cL(\cF^2(E,\si))$ for $b\in\si(\cA)'$.

In a similar way as above, we define the Fock spaces
$\cF^2(E^\si)=\oplus_{n=0}^\infty (E^\si)^{\otimes n}$ and
$\cF^2(E^\si,\pi)=\cF^2(E^\si)\otimes_\pi\cX$. Moreover, with $\cF^\infty(E^\si,\pi)$ we
denote the space of bounded sequences $\{x(n)\}_{n\in\BZ_+}$ with
$x(n)\in(E^\si)^{\otimes n}\otimes_\pi\cX$.

\begin{example}\label{E:5}(Continuation of Example \ref{E:4}){\em $ $\\
{\bf Part A:} As observed in Example \ref{E:3}, $E^{\otimes n}\otimes_\si\cE=\cU^{d_n}$,
which, via the map $\nu_n$ in \eqref{nu}, can be identified with $\oplus_{\al\in\cF_d,|\al|=n}\cU$.
It thus follows that $\cF^2(E,\si)$ can be identified with $\ell^2_\cU(\cF_d)$, and
similarly $\cF^2(E^\si,\pi)$ with $\ell^2_\cX(\cF_d)$. For both spaces the left
$\BC$-action is just scalar multiplication.
Via the same identification, $\cF^\infty(E^\si,\pi)$ can be seen as the subspace of
$\ell^\infty_\cX(\cF_d)$ consisting of elements $\{x(\al)\}_{\al\in\cF_d}$ from
$\ell^\infty_\cX(\cF_d)$ such that $\sup_{N\in\BZ_+}\sum_{|\al|=N}\|x(\al)\|^2<\infty$.\\
{\bf Part B:} Similarly as in part A, $\cF^2(E,\si)$ and $\cF^2(E^\si,\pi)$
can be identified with $\ell^2_\cK(\cF^d)$ and $\ell^2_\cX(\cF_d)$, respectively,
and $\cF^\infty(E^\si,\pi)$ with a subspace of $\ell^\infty_\cX(\cF_d)$. However,
in this case we have a left $\cL(\cK)$-action that is given by entry-wise
multiplication with $T\in\cL(\cK)$ (resp.~$\pi(T)$) on the left.\\
{}From here on we shall identify the spaces $\cF^2(E,\si)$, $\cF^2(E^\si,\pi)$ and
$\cF^\infty(E^\si,\pi)$, in both parts A and B, with the appropriate $\ell^2$ and
$\ell^\infty$ spaces over $\cF_d$.
}\end{example}

With the state operator $A$ and two nonnegative integers $n,m\in\BZ_+$ we associate
the following {\em generalized power}:
\[
(A_n)^m=(I_{(E^\si)^{\otimes m-1}}\otimes A_n)\cdots(I_{E^\si}\otimes A_n)A_n,
\]
which defines an operator in $\cL((E^\si)^{\otimes n}\otimes\cE,(E^\si)^{\otimes n+m}\otimes\cE)$.
Here $A_n^0$ is to be understood as $I_{(E^\si)^{\otimes n}\otimes\cE}$.
Note that $(A_n)^m=I_{(E^\si)^{\otimes n}}\otimes A_1^m=(A^m)_n$. Thus we can unambiguously
write $A_n^m$ instead of $(A_n)^m$ or $(A^m)_n$.

\begin{theorem}\label{T:dissipative}{\em
Let $y=\{y(n)\}_{n\in\BZ_+}$ and $x=\{x(n)\}_{n\in\BZ_+}$ be the output and
state sequences generated by the dissipative system \eqref{SysEqn} with input sequence
$u=\{u(n)\}_{n\in\BZ_+}$ and initial state $x(0)$.  Assume that $u\in\cF^2(E,\si)$.
Then $y\in\cF^2(E,\si)$ , $x\in\cF^\infty(E^\si,\pi)$ and
\[
y=T_\Si u+W_\Si x(0),
\]
with $T_\Si\in\cL(\cF^2(E,\si))$ and $W_\Si\in\cL(\cX,\cF^2(E,\si))$ the left
$\si(\cA)'$-module maps given by
\begin{equation*}
T_\Si\!=\!\!\mat{ccccc}{
D_0&0&0&\cdots\\
C_1A_1^0B_0&D_1&0&\cdots\\
C_2A_1^1B_0&C_2A_2^0B_1&D_2&\cdots\\
\vdots&\vdots&\vdots&\ddots}\
W_\Si\!=\!\!\mat{c}{C_0A_0^0\\C_1A_0^1\\C_2A_0^2\\\vdots}.
\end{equation*}
Moreover, the operator matrix $\mat{cc}{T_\Si&W_\Si}$ is contractive, and thus
$\|y\|^2\leq\|u\|^2+\|x(0)\|^2$. Finally, if, in addition,
$\|A\|<1$, then $x\in\cF^2(E^\si,\pi)$.
}\end{theorem}

The fact that $T_\Si$ is a $\si(\cA)'$-module map means that it commutes with the operators
$I_{\cF^2(E)}\otimes b$ for each $b\in\si(\cA)'$. In addition, it is not hard to show (see
the proof of Theorem 3.5 in \cite{MSSchur}) that $T_\Si$ also commutes with
$I_{\cF^2(E)}\otimes \mu$ for each $\mu\in E^\si$. In order to make sense of
$I_{\cF^2(E)}\otimes \mu$ as an operator on $\cF^2(E,\si)$, because $\mu$ maps $\cE$ into
$E\otimes\cE$, we have to make the identification $\cF^2(E)\otimes E=\cF^2(E)$. These
commutative relations are one way of defining the analytic Toeplitz algebra $\cT_+(E,\si)$
associated with the CR-pair $(E,\si)$ of \cite{MS98}. There are various other ways to
characterize $\cT_+(E,\si)$; originally in \cite{MS98} the space was generated by
a collection of creation operators, while in \cite{BtH1} a Toeplitz structure characterization
is given. The definition used here stems from \cite{MS04}.

\section{Point-evaluation}\label{S:Point}

Given a dissipative system $\Sys$ as in \eqref{SysMat}, there is a
notion of point evaluation giving rise to a $Z$-transform enabling
the passage from the time domain to the frequency domain for this
abstract setting. In the $W^*$-correspondence setting studied here,
points are pairs $(\eta,b)$ with $\eta$ an element of the {\em
generalized disk}
\[
\BD(E^{\si*}):=\{\eta\colon \eta^*\in E^\si,\ \|\eta\|<1\}.
\]
and $b\in\si(\cA)'$. With a sequence $f=\{f(n)\}_{n\in\BZ_+}\in\cF^2(E,\si)$ (so
$f(n)\in E^{\otimes n}\otimes \cE$), we associate a function $f^\wedge$
on $\BD(E^{\si*})\times\si(\cA)'$ with values  in $\cE$ defined by
\[
f^\wedge(\eta,b)=\sum_{n=0}^\infty \eta^n(I_{E^{\otimes n}}\otimes b)f_n\quad
((\eta,b)\in\BD(E^{\si*})\times\si(\cA)').
\]
Here we use $\eta^n$ to indicate the $n^\tu{th}$ {\em generalized power} of $\eta$:
\[
\eta^n=\eta (I_{E} \otimes \eta) \cdots(I_{E^{\otimes n-1}}\otimes \eta)
:E^{\otimes n}\otimes\cE\to\cE,
\]
where for $n=0$ we set $\eta^0=I_\cE$. Note that the generalized powers
$\eta^n$ are well defined because $\eta$ is a left $\cA$-module map; $\eta$
inherits this property from the fact that $\eta^*\in E^\si$ is a left
$\cA$-module map. The sum in the point evaluation is well defined because
$\|\eta^n\|\leq\|\eta\|^n$ for each $n$.

With a element $x=\{x(n)\}_{n\in\BZ_+}\in\cF^\infty(E^\si,\pi)$ we also
associate a function $x^\wedge$ on $\BD(E^{\si *})\times\si(\cA)'$, but
this point evaluation is essentially different from that in $\cF^2(E,\si)$.
Given $\mu\in E^\si$, we define the
operator $L_\mu:\si(\cA)'\to E^\si$ by
\[
L_\mu b=\mu\otimes b=\mu b\quad (b\in\si(\cA)').
\]
Then $L_\mu$ is a right $\si(\cA)'$-module map, and thus, for $n=1,2,3,\ldots$, we can define the
{\em $n^{th}$ generalized power:}
\[
L_\mu^n=(L_\mu\otimes I_{(E^\si)^{\otimes n-1}})\cdots (L_\mu\otimes I_{E^\si})L_\mu,
\]
which is an operator from $\si(\cA)'$ into $(E^\si)^{\otimes n}$; we set $L_\mu^0$ equal
to the identity map on $\si(\cA)'$. We further define $L_{\mu,\pi}=L_\mu\otimes I_\cX$ and
$L_{\mu,\pi}^n=L_\mu^n\otimes I_\cX$. Then one can compute that
$\|L_{\mu,\pi}^n\|\leq \|\mu\|^n$ for each $n\in\BZ_+$.
In particular, for each $\eta\in\BD(E^{\si *})$
we have the operators $L_{\eta*,\pi}$ and $L_{\eta*,\pi}^n$.
Then, with a state sequence $x=\{x(n)\}_{n\in\BZ_+}\in\cF^\infty(E^\si,\pi)$ we associate a
function $x^\wedge$ on $\BD(E^{\si *})\times \si(\cA)'$ with values in $\cX$ that is defined by
\[
x^\wedge(\eta,b)
=\sum_{n=0}^\infty (L_{\eta^*,\pi}^{n})^*(b\otimes I_{(E^\si)^{\otimes n-1}\otimes\cX})x(n),
\]
for each $(\eta,b)\in\BD(E^{\si *})\times \si(\cA)'$.
In a different context, and restricted to the case $x=\{x(n)\}_{n\in\BZ_+}\in\cF^2(E^\si,\pi)$,
this form of point evaluation was first studied in \cite{MSPoisson}; see also \cite{BtH1}.

\begin{example}[Continuation of Example \ref{E:5}]\label{E:6}{\em $ $\\
{\bf Part A:} In this case $\BD(E^{\si*})$ can be identified with the unit ball
$\BB^d$ defined in \eqref{BB} (formally elements of  $\BD(E^{\si*})$ are row vectors)
and $\si(\cA)'=\BC$.
Let $\eta =(z_1,\ldots,z_d)\in\BB^d$. Then the generalized power $\eta^n$ corresponds
to the vector in $\BC^{d_n}$ whose $j^{th}$ entry is equal to
\[
z^\al=z_{i_n}\cdots z_{i_1}\text{ if }\al:=\nu_n(j)=i_n\cdots i_1.
\]
Since $z_1,\ldots,z_d$ are scalars, and thus commute, the entries in $\eta^n$ corresponding to
$\al_1$ and $\al_2$ (via the map $\nu_n$) that are mapped to the same element in $\BZ_+^d$
via the abelianization map $\Abel$ in \eqref{abel} are the same. Then for $u=\{u(\al)\}_{\al\in\cF_d}$
in $\cF^2(E,\si)=\ell^2_\cU(\cF_d)$, the point evaluation in a point $(z,b)\in\BB^d\times\BC$
is given by
\[
\widehat{u}(z,b)=\sum_{\al\in\cF_d}z^\al b u(\al)=b\sum_{n\in\BZ_+^d} z^n\ov{u}(n)
\]
where $\ov{u}(n)=\sum_{\Abel(\al)=n}u(\al)$ for each $n\in\BZ_+^d$. In particular, the part of the
point evaluation that comes from $\si(\cA)'$ plays no real role, and will be left out in
the sequel, i.e., we only consider $b=1\in\BC$ and shall write $\widehat{u}(z)$
instead of $\widehat{u}(z,1)$. Using that $\cF^\infty(E^\si,\pi)$ can be identified with a
subspace of $\ell^\infty_\cX(\cF_d)$, it follows in a similar way that the point-evaluation
for an $x=\{x(\al)\}_{\al\in\cF_d}$ in $\cF^2(E^\si,\pi)$ is defined by:
\[
\widehat{x}(z)=\sum_{n\in\BZ_+^d} z^n\ov{x}(n)\quad (z\in\BB^d),
\]
where $\ov{x}(n)=\sum_{\Abel(\al)=n}x(\al)$ for each $n\in\BZ_+^d$.\\
{\bf Part B:} In the noncommutative setting $\BD(E^{\si*})$ corresponds to
the set of strict row contractions $T=(T_1,\ldots,T_d)$ with $T_j\in\cL(\cK)$
and $\si(\cA)'=\cL(\cK)$. Similar to the case for Part A, the $n^{th}$ generalized power
of $T$ is defined via a noncommutative functional calculus:~the $j^{th}$ entry
of $T^n$ is given by $T^\al$ in \eqref{NonComFunCalc} in case $\al:=\nu_n(j)$. Since the
elements $T_j$ of $T$ do not (necessarily) commute pairwise, we obtain a
noncommutative point evaluation for an element $u=\{u(\al)\}_{\al\in\cF_d}$ in
$\cF^2(E,\si)=\ell^2_\cK(\cF_d)$:
\[
\hat{u}(T,b)=\sum_{\al\in\cF_d} T^\al b u(\al)
\]
where $T$ is a strict row contraction and $b\in\cL(\cK)$. The point evaluation of an
element $x=\{x(\al)\}_{\al\in\cF_d}$ in $\cF^\infty(E^\si,\pi)=\ell^2_\cX(\cF_d)$ goes
through the $*$-representation $\pi$:
\[
\hat{x}(T,b)=\sum_{\al\in\cF_d} \pi(T^\al b) x(\al)
\]
with again $T$ a strict row contraction and $b\in\cL(\cK)$. Now, under the restrictions
of Example \ref{E:4} we have input sequences of the form $u=\{\ov{u}(\al)k_0\}_{\al\in\cF_d}$,
with $\ov{u}(\al)\in\BC$, and state sequences of the form
$x=\{\ov{x}(\al)\otimes k_0\}_{\al\in\cF_d}$, with $\ov{x}(\al)\in\ov{\cX}$, where $k_0$ is
a fixed vector in $\cK$. In that case $\hat{u}(T,b)=\sum_{\al\in\cF_d}\ov{u}(\al)T^\al b k_0$
and $\hat{x}(T,b)=\sum_{\al\in\cF_d}\ov{x}(\al)\otimes (T^\al b k_0)$,
and replacing $k_0$ with $bk_0$ we see that we may restrict the point evaluation to the
case that $b=I_\cK$.
}\end{example}

\section{Frequency-domain equations}\label{S:Frequency domain equations}

With the point evaluations as defined in the previous section we have the following result.

\begin{theorem}\label{T:FreqEqn}{\em
Let $\Sys$ be a dissipative system. Given an input sequence
$u=\{u(n)\}_{n\in\BZ_+}\in\cF^2(E,\si)$ and an initial state $x(0)$,
let $y=\{y(n)\}_{n\in\BZ_+}\in\cF^2(E,\si)$ and
$x=\{x(n)\}_{n\in\BZ_+}\in\cF^\infty(E^\si,\pi)$ be
the output and state sequences given by the system of equations \eqref{SysEqn}.
Then the $Z$-transforms $u^\wedge$, $x^\wedge$ and $y^\wedge$ of $u$, $x$
and $y$ satisfy the frequency domain equations:
\begin{equation}\label{FreqEqn}
\begin{array}{rcl}
x^\wedge(\eta,b)
\!\!\!\!&=&\!\!\!\!x(0)+L_{\eta^*,\pi}^*Ax^\wedge(\eta,b)+L_{\eta^*,\pi}^*Bu^\wedge(\eta,b),\\[.2cm]
y^\wedge(\eta,b)
\!\!\!\!&=&\!\!\!\!Cx^\wedge(\eta,b)+Du^\wedge(\eta,b),
\end{array}
\end{equation}
and we have
\begin{equation}\label{FreqTransfer}
y^\wedge(\eta,b)=F_\Si(\eta)u^\wedge(\eta,b)+W_\Si(\eta)x(0),
\end{equation}
where the transfer function $F_\Si$ and observation function $W_\Si$ are given by
\begin{equation}\label{TransObs}
\begin{array}{c}
F_\Si(\eta)=D+C(I-L_{\eta^*,\pi}^*A)^{-1}L_{\eta^*,\pi}^*B,\\[.2cm]
W_\Si(\eta)=C(I-L_{\eta^*,\pi}^*A)^{-1}.
\end{array}
\end{equation}
}\end{theorem}

\begin{example}[Continuation of Example \ref{E:6}]\label{E:7}{\em $ $\\
{\bf Part A:}   For this example the system matrix has the form \eqref{ExSysmat}.
An element $\eta$ of ${\mathbb D}(E^{\sigma*})$ is identified with a point
$z = (z_1  , \dots, z_d) \in {\mathbb B}^d$, and in that case the operator
$L_{\eta^*}^* A \in \cL(\cX)$ works out to be $z_1 A_1 + \cdots + z_d A_d$.
Similarly $L_{\eta^*}^* B = z_1 B_1 + \cdots z_d B_d$ and \eqref{TransObs}
gives us
\begin{align*}
 F_\Sigma(z) &= D + C (I - z_1A_1 - \cdots - z_d A_d)^{-1} \times \\
 & \quad  \times (z_1B_1 + \cdots + z_dB_d) \\
 W_\Sigma(z) &= C(I - z_1A_1 - \cdots - \cdots z_d A_d)^{-1}
\end{align*}
for the transfer function and the observation function, in agreement
with \eqref{CFMtrans-obs}, and
\eqref{FreqTransfer} reduces to \eqref{Cinput-output}.

{\bf Part B:}  For this example,  $\sigma(\cA)' = \cL(\cK)$ is
nontrivial, $\pi \colon \si(\cA)' = \cL(\cK) \to \cL(\cX)$ is a nontrivial
representation, and ${\mathbb D}(E^{\sigma *})$ is
identified with the set of strict row contractions $T = \begin{bmatrix} T_1 &
\cdots & T_d \end{bmatrix} \in \cL(\cK^d, \cK)$.  Nevertheless $E^\sigma \otimes_\pi \cX$ is
again identified with $\cX^d$ and the system matrix has the form $\eqref{ExSysmat}$.
Following the assumptions made in Example \ref{E:4}, i.e., $\cX=\ov{\cX}\otimes\cK$
and the system operators are of the form \eqref{inertwine2},
if $\eta\in {\mathbb D}(E^{\sigma *})$ is identified with the strict row
contraction $T = \begin{bmatrix} T_1 & \cdots & T_d \end{bmatrix}$,
then it works out that the operators $L_{\eta^*}^* A  \in \cL(\cX)$
and  $L_{\eta^*}^* B  \in \cL(\cE, \cX)$ can be identified with
\begin{align*} & L_{\eta^*}^*A =  \ov{A}_1 \otimes T_1 + \cdots  +  \ov{A}_d \otimes T_d,  \\
&  L_{\eta^*}^* B  = \ov{B}_1 \otimes T_1 + \cdots  + \ov{B}_d \otimes T_d
 \end{align*}
whence the transfer function and observation function assume the form
\begin{align*}
F_\Sigma(T) & = D + \ov{C}\otimes I_\cK
(I - \ov{A}_1 \otimes T_1 - \cdots - \ov{A}_d \otimes T_d)^{-1} \times \\
 & \qquad \qquad  \times (\ov{B}_1 \otimes T_1 + \cdots \ov{B}_d \otimes T_d), \\
 W_\Sigma(T) & = \ov{C}\otimes I_\cK
 (I - \ov{A}_1 \otimes T_1 - \cdots - \ov{A}_d \otimes T_d)^{-1},
 \end{align*}
 in agreement with \eqref{NCtrans-obs}, and, moreover, \eqref{FreqTransfer} collapses to \eqref{ii-o}.
}\end{example}

\section{Other examples and other types of generalized systems}  \label{S:other}

Our main purpose is to illustrate the general Muhly-Solel framework
of tensor algebras, generalized Schur class, and Hardy algebras
arising from a CR pair for two specific cases, namely: (1)
commutative Fornasini-Marchesini input/state/output systems and
their associated transfer functions as holomorphic functions on the
unit ball, and (2) noncommutative Fornasini-Marchesini systems with
associated transfer functions defined on the noncommutative unit
ball.  We mention that the Muhly-Solel formalism has at least one
other special case which is a familiar setting in system theory,
specifically, time-varying systems with the Muhly-Solel Schur-class
function identifiable as the input-output map for a discrete-time
linear time-varying system.  Details on this can be found in
\cite{BBFtH08}.

Still other types of systems  and associated transfer functions have
been introduced and studied recently.  Specifically there is a
system theory over hyper-holomorphic functions (with quaternions
playing the role of the complex numbers) and a system theory for
stochastic systems where the system parameters are themselves random
variables;  in these instances the familiar algebra centering around
pointwise multiplication of functions is replaced by something more
exotic:  Cauchy-Kovalevskaya product of hyperholomorphic functions
or Wick product of stochastic variables.   An introductory survey to
these matters can be found in \cite{Alpay}.


\end{document}